\documentclass[12pt]{article}
\usepackage{amsfonts}
\begin{document}
\def\R{\mathbb{R}}
\def\C{\mathbb{C}}
\def\Z{\mathbb{Z}}
\def\N{\mathbb{N}}
\def\Q{\mathbb{Q}}
\def\D{\mathbb{D}}
\def\E{\mathbb{E}}
\def\Sp{{\mathbb{S}}}
\def\T{\mathbb{T}}
\def\hb{\hfil \break}
\def\ni{\noindent}
\def\i{\indent}
\def\a{\alpha}
\def\b{\beta}
\def\e{\epsilon}
\def\d{\delta}
\def\De{\Delta}
\def\g{\gamma}
\def\qq{\qquad}
\def\L{\Lambda}
\def\G{\Gamma}
\def\F{\cal F}
\def\K{\cal K}
\def\A{\cal A}
\def\B{\cal B}
\def\M{\cal M}
\def\P{\cal P}
\def\Om{\Omega}
\def\om{\omega}
\def\s{\sigma}
\def\th{\theta}
\def\Th{\Theta}
\def\z{\zeta}
\def\p{\phi}
\def\m{\mu}
\def\n{\nu}
\def\la{\lambda}
\def\Si{\Sigma}
\def\q{\quad}
\def\qq{\qquad}
\def\half{\frac{1}{2}}
\def\hb{\hfil \break}
\def\half{\frac{1}{2}}
\def\pa{\partial}
\def\r{\rho}
\def\Spd{{{\Sp}}^d}
\begin{center}
{\bf TAUBERIAN KOREVAAR}
% {\bf SOME ASPECTS OF JAAP KOREVAAR'S WORK ON TAUBERIAN THEORY}
% {\bf TAUBERIAN THEOREMS, KOREVAAR AND ME} 
\end{center}
\begin{center}
{\bf N. H. BINGHAM}
\end{center}

\ni {\bf Abstract}  We focus on the Tauberian work for which Jaap Korevaar is best known, together with its connections with probability theory.  We begin (\S 1) with a brief sketch of the field up to Beurling's work.  We follow with three sections on Beurling aspects: Beurling slow variation (\S 2); the Beurling Tauberian theorem for which it was developed (\S 3); Riesz means and Beurling moving averages (\S 4).  We then give three applications from probability theory: extremes (\S 5), laws of large numbers (\S 6), and large deviations (\S 7).  We turn briefly to other areas of Korevaar's work in \S 8.  We close with a personal postscript (whence our title). \\

\ni {\it Keywords}  Tauberian theorem, Beurling slow variation, Beurling Tauberian theorem, Beurling moving average, Riesz mean, large deviations \\
\ni {\it MSC Classification} 01A70, 40E05    \\    

\ni {\bf 1. Tauberian theorems: from Tauber to Beurling} \\
% through Hardy, Littlewood, Karamata and Wiener
\i The origins of Tauberian theory can be traced to Abel's continuity theorem of 1826, which we teach to students.  In 1828, Abel described divergent series as the invention of the devil, and denounced any use of them.  But, it was found that some divergent series could be made convergent by suitable smoothing or averaging, leading to what we now call summability theory(see e.g. [ZelB1970]).  Then in 1897 Tauber [Tau1897] gave the first of what Hardy and Littlewood later named a {\it Tauberian theorem}: a `corrected converse' leading from convergence of a smoothed sequence to that of the original sequence, under a suitable condition, a {\it Tauberian condition}.  In his Preface to Hardy's {\sl Divergent series} [Har1949], Littlewood comments of the resulting theory that `in the early years of the century the subject, while in no way mystical or unrigorous, {\it was} regarded as sensational, and about the present title, now colourless, there hung an aroma of paradox and audacity'. \\
\i At first, attention focused on special methods of summation, such as the Abel, Ces\`aro, Euler and Borel [Har1949, I-XI].  One culminating result was the Hardy-Littlewood-Karamata (Tauberian) theorem of 1930 for Laplace transforms [Kor2004, I].  The whole area changed with the introduction of Wiener's general Tauberian theory [Wie1932]; [Har1949, XII], [Kor2004, II].  One of Korevaar's early contributions here was his strikingly short proof of Wiener's Tauberian theorem by Schwartz distributions (generalised functions) [Kor1965]. \\
\i The {\it Borel} (exponential) {\it method} $B$ of summability is defined by writing 
$$
s_n \to s \qq (n \to \infty) \qq (B)
$$
(handling series $\sum a_n$ and sequences $s_n$ together by $s_n := \sum_0^n a_k$) for 
$$
e^{-x} \sum_0^{\infty} s_n x^n/n! \to s \qq (x \to \infty).
$$
This method is useful for analytic continuation of power series [Har1949, VIII, IX], [Kor2004, VI] and in number theory [Ten1980].  It is also useful in probability theory, since the weights $e^{-x} x^n/n!$ are those of a Poisson law $P(x)$ with parameter (mean, variance) $x$.  It is intimately linked with the {\it Valiron method} $V_{\half}$ [BinT1986] (or just $V$ here for convenience, defined by writing      
$$
s_n \to s \qq (n \to \infty) \qq (V)
$$
for
$$
\frac{1}{\sqrt{2 \pi x}} \sum_0^{\infty} s_n \exp \{ - \half (x - n)^2/x \} \to s \qq (x \to \infty)
$$
[HarL1916].  Such a link is probabilistically indicated by the central limit theorem, as for large $x$ the Poisson law approximates a normal, $P(x) \sim N(x,x)$. \\

\ni {\bf 2.  Beurling slow variation} \\ 
\i For reasons we turn to below regarding his Tauberian theorem, Beurling (in unpublished lectures of 1957, published by others in 1972 [Kor2004, \S IV.11]) introduced what is now called {\it Beurling slow variation}.  This is a variant on Karamata slow variation (for which see e.g. [BinGT1987], `BGT'), and as there can be formulated for measurable functions (as was done by Korevaar et al. [KorvA1949]), or topologically for functions with the Baire property [Mat1964] (in brief, Baire functions), and indeed more generally. \\  
\i A positive Baire or measurable (Baire/measurable below) function ${\phi}$ on ${\R}_+$ is called {\it Beurling slowly varying} if $\phi(x) = o(x)$ as $x \to \infty$ and
$$
\phi(x + t \phi(x))/\phi(x) 
\to 1 \qq (x \to \infty) \q \forall t.                      \eqno(BSV)
$$
Writing
$$
x \ {\circ}_{\phi} \ t := x + t \phi(x),
$$
this is
$$
\phi(x \ {\circ}_{\phi} \ t)/\phi(x) \to 1 \q \forall t.
$$
In recent papers (see e.g. [BinO14], [BinO16], [BinO20a]) the author and Ostaszewski take the view that it is (Baire) {\it category} rather than {\it measure} which is the principal case, the reverse of the chronological order, whence our ordering above.   \\
\i By Bloom's theorem [Blo1976], for $\phi$ Beurling slowly varying ($\phi \in BSV$) and  {\it continuous}, the convergence in $(BSV)$ is locally uniform (uniform on compact $t$-sets in ${\R}_+$).  Continuity is weakened to the {\it Darboux} (or intermediate-value) property in [BinO2014].  There, a number of properties are shown to be sufficient, including monotonicity, but whether continuity may simply be dropped here remains open. \\
\i When the convergence in $(BSV)$ is is locally uniform, $\phi$ is called {\it self-neglecting}, $\phi \in SN$. \\
\i For $\phi \in SN$, a Baire/measurable $f$ is {\it Beurling regularly varying} [BinO2014] with {\it limit} $g$ if
$$
f(x + t \phi(x))/f(x) \to g(t) \qq (x \to \infty) 
\qq \forall \ t.                                          \eqno(BRV)
$$
\i The condition $\phi(x) = o(x)$ may be usefully weakened to $\phi(x) = O(x)$ [Ost2015].  For $\phi(x) = O(x)$, call $\phi > 0$ {\it self-equivarying} with {\it limit} $\la$ if 
$$
\phi(x + t \phi(x))/\phi(x) \to \la(t) \q 
\hbox{locally uniformly}.                                   \eqno(SE)
$$
The limit $\la$ then satisfies the {Go{\l}ab-Schinzel functional equation}
$$
\la(x) \la(y) = \la(x + y \la(x)) \qq \forall \ x, y.      \eqno(GS)
$$
It is continuous, and of the form
$$
\la(t) = 1 + at
$$
for some $a \geq 0$.  To within re-scaling, there are only two cases: the `small-order limit' $\la(t) \equiv 1$, and `large-order limit' $\la(t) = 1 + t$ [Ost2015, Th. 0]. \\ 
\i We note that as $\phi(x) = O(x)$, $\int_1^x dt/\phi(t)$ diverges (logarithmically) for large $x$.  Its unboundedness will be important in \S 4 in connection with Riesz means. \\   
  
\ni {\bf 3. Beurling's Tauberian theorem} \\
\i The Borel method [Bor1899], while very useful and interesting, is notably more difficult to handle than most of the other classical summability methods in common use (Abel, Ces\`aro, Euler etc.), witness its extended treatment in [Har1949] and [Kor2004].  Beurling's Tauberian theorem, to which we now turn, was developed to extend the Wiener Tauberian theorem conveniently to cover the Borel method. \\
\i Recall that $K \in L_1(\R)$ satisfies the {\it Wiener condition}, or is a {\it Wiener kernel}, if its Fourier transform $\hat K(t)$ has no real zeros $t$.  For $H$ bounded and $F \in L_1(\R)$, define the {\it Beurling convolution} (w.r.t. $\phi$ as in (BSV)) by
$$
F *_{\phi} H (x) := 
\int_{\R} F \Bigl(\frac{x - u}{\phi(x)} \Bigr) H(u) \frac{du}{\phi(x)}
= \int F(-y) H(x \ {\circ}_{\phi} \ t) dt.
$$

\ni {\bf Theorem B (Beurling's Tauberian theorem, 1957)}.  For $K$ a Wiener kernel, $\phi$ Beurling slowly varying: if $H$ is bounded and
$$
K *_{\phi} H (x) \to c \int K(y) dy \qq (x \to \infty),
$$
then for all $F \in L_1(\R)$,
$$
F *_{\phi} H (x) \to c \int F(y) dy \qq (x \to \infty).
$$

\i The case $\phi \equiv 1$ gives the Wiener Tauberian theorem, Theorem W say ([Wie1932]; BGT \S 4.8, [Kor2004,II]): \\

\ni {\bf Wiener's Tauberian theorem, 1932}.  Theorem B holds with $\phi$-convolution $*_{\phi}$ replaced by ordinary convolution $*$. \\

\i Recall that Theorem W is proved from the {\it Wiener approximation theorem} (or {\it closure} theorem): $K$ is a Wiener kernel if and only if linear combinations of its translates $K(a \ + \ .)$ are dense in $L_1(\R)$ (see e.g. [Kor2004, V.3]).  Korevaar [Kor2004, IV.11] gives a very short proof of Theorem B by reducing it to the Wiener approximation theorem. \\
\i `Wiener's second Tauberian theorem', involving Lebesgue-Stieltes integrals [Kor2004, II.13], also has a Beurling counterpart [BinO16, \S 4]. \\ 

\ni {\bf 4.  Riesz means and Beurling moving averages} \\
\i First, recall the {\it Riesz (typical) means} $R(\la,1)$ of order 1, or simply $R(\la)$ here as we shall only need order 1, where $\la := ({\la}_n)$ is a real sequence ${\la}_n \uparrow \infty$: we write
$$
s_n \to s \qq R(\la)
$$
for
$$
\frac{1}{x} \int_0^x \Bigl( \sum_{{\la}_n \leq y} a_n \Bigr) dy \to s \qq (x \to \infty).                                      \eqno(R(\la))
$$
\i The first consistency theorem for Riesz means [ChaM1952, I] relates Riesz means of different orders, so need not detain us as we use only order 1 here.  The {\it second consistency theorem for Riesz means} (Riesz in 1909 and 1916, Hardy in 1910 and 1915; see [ChaM1952, II]), important for us, tells us that two Riesz means $R({\la})$, $R({\mu})$ whose {\it logarithms are comparable}, 
$$
0 < \liminf \log {\la}_n/ \log {\mu}_n \leq \limsup \log {\la}_n/\log {\mu}_n < \infty,                                  \eqno(log \lambda)
$$
are equivalent.  Thus the {\it order of magnitude of the logarithm} is all that matters. \\   
\i Since a $\phi \in SN$ is Beurling slowly varying, in
$$
\frac{1}{t \phi(x)} \sum_{x < n \leq x + t \phi(x)} s_n \to s
\qq (x \to \infty) \qq \forall \ t                             \eqno(BMA)
$$
the left-hand side is called a {\it Beurling moving average}.  Then ([BinGo1988]; cf. [BinT1986]), writing
$$
\la(x) := \exp \{ \int_1^x dt/\phi(t) \},
$$ 
$(R(\la))$ and $(BMA)$ are equivalent. \\
\i As noted in \S 2, all $\int_1^x dt/\phi(t) \uparrow \infty$ as $x \to \infty$.  In view of the second consistency theorem, any choice of sequence $\la := ({\la}_n), {\la}_n \uparrow \infty$ agreeing with the function $\la(x)$ above at all points $x = {\la}_n$ give equivalent Riesz means, and so we are free to choose any and then abbreviate $\la(.)$ to $\la$ without ambiguity.  More is true: agreement may be weakened to approximation within reasonable limits, subject only to 
$(log \lambda)$.  \\
\i The title of [Bin2019] coincides with that of this section, and we can refer to it (or its ArXiv version) for more detail than we give here.  It has been known since [Bin1981] that Beurling moving averages are special cases of Riesz means.  In particular, [Bin2019 Th. 8.3] gives a transparent proof of the general case of results in [BinT1986 Th. 3] and [BinGo1983, Th. 5]: that convergence as in $(BMA)$ is equivalent to (ordinary) convergence of an $o(1)$-{\it perturbation} with {\it rate} $1/\phi(x)$.  As noted in [Bin2019], this result, though transparent when seen as a representation theorem for a form (Beurling) of regular variation, has distinguished antecedents, going back to Hardy in 1904 [Har1949, Th. 149]; see [Bin2019, \S 2] for its full history. \\      

\ni {\bf 5.  Extremes} \\
\i The mathematics (probability and statistics) of extremes is a vast and topical area, recently surveyed in [BinO2021], so we can be brief here, referring there for detail and references.  The subject essentially dates from Fisher and Tippett in 1928.  There, they obtained the three classical limit laws of maxima $M_n := \max \{ X_1, \cdots, X_n \}$ of independent and identically distributed (iid) random variables, which are to within type (location and scale) the Fr\'echet (heavy-tailed, ${\Phi}_{\a}$, $\a > 0$), Gumbel (light-tailed, $\Lambda$) and Weibull (bounded tail, ${\Psi}_{\a}$, $\a > 0$).  The corresponding domains of attraction (laws $F$ sampling from which can give such a limit) $D({\Phi}_{\a})$, $D(\Lambda)$, $D({\Psi}_{\a})$ are simple to describe in the Fr\'echet case (in terms of Karamata regular variation, `at infinity') and the Weibull case (regular variation `at the finite end-point').  The Gumbel domain of attraction is more complicated, and best described nowadays in the language of Beurling regular variation (below).  We note that Fisher and Tippett (who gave no references and merely sketches of proofs) noticed and studied in detail that for $F = \Phi$ the (standard) normal, the convergence is {\it extremely slow} -- so slow that for numerical purposes it is better to avoid the `ultimate' approximation (to the Gumbel) and use instead the `penultimate' approximation (to the Fr\'echet).  This was both astute theoretically and impressive numerically in the pre-computer days of desk machines (the penultimate approximation has been developed more recently by R. L. Smith and J. P. Cohen). \\
\i For simplicity, we work to within type (centre and scale, to have mean 0 and variance 1).  We may then combine the three cases above into a single one-parameter family, the {\it extreme-value distributions} (EVD) $G_{\a}, \ \a \in \R$, where $\a > 0$ for Fr\'echet, $\a = 0$ for Gumbel, $\a < 0$ for Weibull, using the `L'Hospital convention': as $(1 + x/n)^n \to e^x$ as $n \to \infty$, we interpret the $\a = 0$ case of $(1 + \a x)^{1/\a}$ as $e^x$.  This gives
$$
G_{\a}(x) := \exp (- g_{\a}(x)), \qquad 
g_{\a}(x) := [ 1 + \a x]_+^{-1/{\a}}.                                  \eqno(EVD)
$$ 
Here the parameter $\a \in \R$ is called the {\it extreme-value index (EVI)} or {\it extremal index}.  The upper end-point $x_+$ of $F$ is $\infty$ for $\a \geq 0$ (with a power tail for $\a > 0$ and an exponential tail for $\a = 0$); for $\a < 0$ $x_+ = -1/\a$, with a power tail to the left of $x_+$. \\
\i The Gumbel domain of attraction, due to de Haan in 1970-71 (see e.g. BGT, Th. 8.13.4; cf. BGT, Ch. 3, De Haan theory) is given by
$F \in D(\Lambda)$ iff 
$$
\overline{F}(t + x a(t))/\overline{F}(t) 
\to g_0(x) := e^{-x} \qquad (t \to \infty),               \eqno(\ast)           
$$
for some {\it auxiliary function} $a > 0$, $a \in SN$, which may be taken [EmbK, (3.34)] as
$$
a(t) := \int_t^{x_+} \overline{F}(u) du/\overline{F}(t) 
\quad (t < x_+),                                          \eqno(aux)       
$$
and satisfies (in the usual case, $x_+ = \infty$)
$$
a(t + x a(t))/a(t) \to 1 \quad (t \to \infty), \qq \hbox{locally uniformly}.                                                \eqno(Beu)
$$
\i The three domain-of-attraction conditions may be unified (again using the L'Hospital convention): $F \in D(G_{\a})$ iff
$$
\overline{F}(t + x a(t))/\overline{F}(t)
\to g_{\a}(x) := (1 + \a x)_+^{-1/\a} \quad (t \to \infty)                   \eqno(\ast \ast)
$$
for some auxiliary function $a$, and then
$$
a(t + x a(t))/a(t) \to 1 + \a x \quad (t \to \infty),           \eqno(\a Beu)
$$
extending the $\a = 0$ case $(Beu)$ above  (see e.g. [BieG, \S 2.6]). \\
\i As explained in [BinO2021], there is no essential loss in assuming (following von Mises) that the density $f$ of the law $F$ exists, in which case the inverse hazard function $i = 1/h$ of the hazard function $h$,
$$
i(t) := \int_t^{\infty} f(u) du/f(x),
$$
may be used as auxiliary function in place of $a$ if preferred.  So may     the mean excess function
$$
e(t) := \E[X - t | X > t].
$$
\i For an application of Theorem B to records, see El Arrouchi [ElA2017]. \\     

\ni {\bf 6.  Laws of large numbers} \\
\i For $X, X_1, \cdots, X_n, \cdots$ iid, we recall Kolmogorov's classic strong law of large numbers (his {\sl Grundbegriffe}, 1933): \\  
$$
\E[ \ |X| \ ] < \infty \ \& \ \E[X] = \mu \ \ \Leftrightarrow \ \  
\frac{1}{n} \sum_1^n X_k \to \mu \q (n \to \infty) \q a.s.
$$
\i The summability method used here is the Ces\`aro $C_1$.  This is embedded in a one-parameter family $\{ C_{\a}: \a > 0 \}$ (see e.g. [Har1949, V, VI]; $\a > -1$ is possible), ordered by inclusion: increasing $\a$ sums more series but gives a weaker conclusion.  Here
$$
s_n \to s \qq (C_{\a})
$$
means
$$
\frac{1}{A_n^{\a}} \sum_0^n A_{n-k}^{\a - 1} s_k \to s
$$
(using the full range $\a > -1$ here), where
$$
A_n^{\a} := (\a + 1) \cdots (\a + n)/n! \sim n^{\a}/\Gamma(1 + \a) \qq (n \to \infty).
$$  
\i The key role of Kolmogorov's strong law is reflected in a striking {\it discontinuity} of the behaviour of $C_{\a}$ here across $\a = 1$ [Bin1989, Th. 1]: \\

\ni {\bf Theorem K}. (i) For $0 < \a \leq 1$,
$$
X_n \to \mu \ \ a.s. \ \ (C_{\a}) \ \ \Leftrightarrow \ \ 
\E[ \ |X|^{1/\a} \ ] < \infty \ \ \& \ \ \E[X] = \mu;
$$
(ii) For $\a  \geq 1$,
$$
X_n \to \mu \ \ a.s. \ \ (C_{\a}) \ \ \Leftrightarrow \ \ 
\E[ \ |X| \ ] < \infty \ \ \& \ \ \E[X] = \mu.
$$

\ni {\it Proof}. (i).  The case $\a = 1$ is Kolmogorov's result (so, Theorem K in his honour). \\
\i For $\a \in (\half, 1)$ ($p := 1/\a \in (1,2)$), this is Theorem 3 of Lorentz [Lor1955]. \\
\i For $\a \in (0, \half)$ ($p := 1/\a > 2$), this is the special case $c_n = A_n^{\a - 1}$ of Theorem 1 of Chow and Lai [ChoL1973]: in the zero-mean case, $X \in L_{1/\a}$ if and only if
$$
n^{-\a} \sum_0^n c_{n-k} X_k \to 0 \qq a.s.
$$
for some (all) $c = (c_n) \in {\ell}_2$ (as $\a - 1 < - \half$ here). \\
\i For $\a = \half$, the result is due to D\'eniel and Derriennic [DenD]. \\
(ii) This is contained in Lai's Theorem L below [Lai1974], that in this iid setting for $\a \geq 1$ the Ces\`aro methods $C_{\a}$ are all equivalent, to each other and to the Abel method A. \hfil $\square$ \break
\\
\ni {\bf Theorem L (Lai, 1974)}.  The following are equivalent: 
$$
\E[ \ |X| \ ] < \infty, \q \E[X] = \mu,
$$
$$
X_n \to \mu \q a.s. \q (C_{\a}) \q \hbox{for some (all)} \ \a \geq 1;
$$
$$
X_n \to \mu \q a.s. \q (A).
$$
\\
\i For the Euler methods $E_p, \ p \in (0,1)$ and the Borel method $B$ [Har1949, VIII, IX], one has: \\
\\
\ni {\bf Theorem C (Chow, 1973)}.  The following are equivalent:
$$
var \ X < \infty, \q \E[X] = \mu,
$$
$$
X_n \to \mu \q a.s. \q (E_p) \q \hbox{for some (all)} \ p \in (0,1);
$$
$$
X_n \to \mu \q a.s. \q (B);
$$
$$
X_n \to \mu \q a.s. \q  (R(e^{\sqrt{n}}));
$$
$$
X_n \to \mu \q a.s. \q  (V_{\half}).
$$ 
\\
\ni {\it Proof}.  The first four statements are in [Cho1973], with the fourth in Beurling moving average (`delayed sum') language, equivalent as above to the Riesz language used here; the fifth is in [Bin1984a, Th. 3]. \hfil $\square$ \break
      
\i Thus the Ces\`aro-Abel family of summability methods corresponds in this setting to means and $L_1$, the Euler-Borel family of methods to variances and $L_2$.  Similarly for the Riesz and Valiron methods (above), the {\it circle methods} (Kreisverfahren: [Har1949, IX], Meyer-K\"onig [MeyK1949]; [Bin1984b]), and the {\it random-walk} methods [Bin1984c].  These are matrix methods $A = (a_{nk})$, whose elements give the distribution of integer-valued random walks $S_n = \sum_1^n X_k$:
$$
a_{nk} = P(S_n = k).
$$
When the step-length law is Poisson, this gives the {\it discrete Borel method}.  Interestingly, whereas the Borel method has a pure gap (high-indices) theorem (Gaier, [Gai1953]; Tur\'an [Tur1984]), the discrete Borel does not (Meyer-K\"onig and Zeller [MeyKZ1960]). \\     
\i For general $p > 1$, one has [BinT1986, Th.5]: \\

\ni {\bf Theorem BT}.  For $p > 1$, the following are equivalent:
$$
\E[ \ | X |^p \ ] < \infty, \qq \E[X] = \mu;
$$
$$
X_n \to \mu \q a.s. \q  (R(\exp (n^{1 - (1/p)}));
$$
$$
X_n \to \mu \q a.s. \q  (V_{1/p}).
$$
\\
\i The results above assume at least as much integrability as existence of the mean.  But one can also assume less [BinGa2015]; here one can no longer centre at means as means do not exist, but must use alternative centring.  Recall the {\it Lambert} $W$-{\it function}, the solution to the functional equation
$$
z = W(z) \exp \{ W(z) \},                                  \eqno(W)
$$
and the {\it logarithmic} summability method $(\ell)$, defined by writing
$$
s_n \to s \qq (\ell)
$$
for 
$$
\frac{1}{\log n} \sum_0^n s_k/(k + 1) \to s \qq (n \to \infty)
$$
(the method $(\ell)$ is equivalent to the Riesz mean $R(\log (n+1))$ [Har1949, Th. 37).  Note that the limit in the Riesz mean is continuous, while that in the logarithmic method is discrete.  This contrast is explored in detail in [BinGa2015] (where other logarithmic methods are also treated), and later in [BinO2020b].  {\it Discontinuous} Riesz means, together with Voronoi means and {\it non-regular} (not necessarily convergence-preserving) summability methods, are considered in [BinGa2017].  \\

\ni {\bf Theorem BG}. For $m_k := \E [X_k I(|X_k| \leq (k+1) \log (k+1)]$, the following are equivalent: \\
(i) $\E[\exp \{ W(|X|) \}] < \infty$, i.e. 
$\E \Bigl[ \frac{|X|}{\log |X| \wedge 1} \Bigr] < \infty$, i.e. 
$\E \Bigl[ \frac{|X|}{1 + {\log}_+ |X|} \Bigr] < \infty$; \\
(ii) $X_n/(n \log n) \to 0$ a.s.; \\
(iii) $X_n - m_n \to 0$ a.s. \\

\i There are numerous other equivalences; applications include the Almost-Sure Central Limit Theorem and the Prime Number Theorem. \\
\i Far-reaching generalizations of results of this type are given in [BinGa2017]. \\ 
 
\ni {\bf 7.  Large deviations} \\
\i We turn now to more detailed results, in which the law $F$, equivalently its characteristic function (Fourier-Stieltjes transform) $\phi(t) := \E[e^{itX}] = \int e^{itx} dF(x)$, has progressively better behaviour.  More moments $m_n := \E [X^k]$ may exist; all moments may exist; the moment sequence may uniquely determine the law $F$; $\phi$ may be analytic (in $t$, now complex, in a neighbourhood of the origin, or a strip in the complex $t$-plane -- Cram\'er's condition); $\phi$ may be entire.  We can say things about events further from typical behaviour the more we progress up this hierarchy,   hence the name large deviations for the area. \\
\i We mention first the Berry-Esseen theorem (see e.g. Petrov [Pet1975, V.2], giving a uniform (in $x$) bound between the law of the sum $S_n := \sum_1^n X_k$ when centred and scaled and the standard normal (Gaussian) $\Phi = N(0,1)$, in terms of the third moment of $F$.  Non-uniform versions are valuable for large $|x|$ [IbrL1971, \S 3.6], [Pet1975, V.4]. \\
\i Expansions related to the central limit theorem when moments beyond the second exist are called Edgeworth expansions (from his work of 1907); see Hall [Hal1992] (who traces them back to Chebyshev in 1890 -- very appropriately, as the area has been dominated by the Russian school; cf. [Bin2021]).  For detailed accounts, see e.g. Feller [Fel1971, XVI], [Pet1975, VI]. \\
\i Returning to the Borel method $B$ of \S 1: Hardy [Har1949, (9.1.8)] gives a detailed approximation of the ratio of Borel weights to Valiron.  In probability language,  for $X \sim P(\la)$, he estimates this in the range $\{ |X - \E[X]| = O({\la}^{\zeta}) \}$, for 
$$
\frac{1}{2} < \zeta < \frac{2}{3}
$$
(Korevaar [Kor2004, 292]) -- the `signature of large deviations' (cf. [IbrL1971, Ch. 12]). \\
\i The Tauberian theorems of exponential type in BGT, \S 4.12 are of large-deviations type; cf. [Kor2004, 208, 292], [Bin2008]. \\            
\i When $\phi(t)$ is entire, one can work instead with the moment-generating function $M(t) = \phi(-it) = \int e^{tx} dF(x)$, also entire, and its logarithm, the cumulant-generatng function $K(t) := \log M(t)$, which is convex.  So one may take its Fenchel dual (or Legendre transform), $K^*$.  By Cram\'er's theorem [Cra1938], writing $F_n$ for the $n$-fold convolution of $F$ (law of $S_n = \sum_1^n X_k$), for measurable sets $A$ with interior $A^{\circ}$ and closure $\bar A$,
\begin{eqnarray*}
- \inf \{ K^*(x): x \in A^{\circ} \} 
&\leq& \liminf \frac{1}{n} \log (F_n(A)) \\
&\leq& \limsup \frac{1}{n} \log (F_n(A)) \\
&\leq& - \inf \{ K^*(x): x \in \bar A \}.
\end{eqnarray*}
One says that the $F_n$ satisfy a {\it large-deviation principle} with {\it rate-function} $K^*$.  An extensive theory of large deviations has developed from this; see e.g. [DeuS1989], [DupE1997], [DemZ1998]. \\ 
\i Estimating the rate of decay of exponentially small probabilities, or of the rate of occurrence of extremely rare events, is important in physics -- e.g., for the half-life of radioactive elements, which determine how long radioactive contamination persists. \\

\ni {\bf 8.  Other areas} \\
\i There is a great deal more to Korevaar the mathematician than being an (or the) expert in Tauberian theorems.  In preparing this centenary tribute, I came across parts of his extensive and impressive corpus that I was not (or was no longer) familiar with.  I mention a few. \\
1. {\it Electrostatics and potential theory}.  Nature likes to arrange herself so to minimise energy (cf. soldiers going from standing at attention to standing at ease, then standing easy).  Korevaar studies how discrete charges approximate the continuous situation (Fekete points). \\
2. {\it Gap (lacunary) series}.  This has been an ongoing interest of Korevaar's.  For other work, we refer to Levinson's book [Lev1940], Mary Weiss's paper [Wei1959], and for probabilistic aspects, Hawkes [Haw1980]. \\ 
3.  {\it Pansions}.  These are `expansions' in Hermite polynomials arising via the Fourier transform, but which are not strictly expansions in any natural sense.  Korevaar [Kor1959] thus calls them pansions; Kahane, in his sympathetic review (MR0104975) in French, correspondingly abbreviates d\'eveloppement to `veloppement'.  (I heard a course of lectures on pansions at the St. Andrews Colloquium in 1968 by de Bruijn, and had forgotten they stem from Korevaar; it was nice to meet them again.) \\
4. {\it Obituaries etc.}  His writings on N. G. (Dick) de Bruijn (1918-2012) [Kor2013] and J. G. van der Corput (1890-1975) [Kor2015] stand out, as do his memories of his teacher H. D. Kloosterman (1900-1968) [Kor2013].  (Adam Ostaszewski and I were happy to contribute to the de Bruijn memorial issue of {\sl Indagationes} [BinO2013].) \\

\ni {\bf Postscript: On Jaap Korevaar} \\    
\i I close by reminiscing briefly about my dealings with Jaap Korevaar over the years.  I fell in love with Tauberian theory in general, and the Wiener and Karamata theories in particular, as a research student in the late 60s, and began to publish on Tauberian theory in the 70s.  I attended the 1979 LMS Durham Symposium on Aspects of Contemporary Complex Analysis, where he gave one of the invited talks.  Early on, when we were wearing our name-badges, I was approached by a smiling, very well-preserved man (whom I now know to have been 56, to my 34), who put his hand out saying ``Tauberian Bingham".  I shook it, saying ``Tauberian Korevaar".  We both roared with laughter, and have been firm friends ever since.  We followed each other's work, and met intermittently at conferences.  In the early years of this century, I had the great good fortune to be consulted by him while he was preparing his splendid magnum opus, his Tauberian book of 2004.  I have his inscribed copy as a treasured possession.  I also had the pleasure of speaking at his 80th birthday conference (my talk, 31 January 2003, `Tauberian theorems, Korevaar and me'), and of meeting his charming wife Pia (Pfluger, daughter of the complex analyst Albert Pfluger), alas, no longer with us. \\
\i Dr Samuel Johnson famously said that every man should strive to be an ornament to his profession.  Jaap Korevaar has most certainly done that, and is an example to us all. \\
\i I thank the editors for their kind invitation to contribute to this centenary volume in his honour.  It gives me great pleasure to do so.\\    

\begin{center}
{\bf REFERENCES}
\end{center}

\ni [BieG2004] Bierlant, J., Goegebeur, Y.,  Segers, J. and Teugels, J. L., {\sl Statistics of extremes: Theory and applications}.  Wiley, 2004. \\
\ni [Bin1981] Bingham, N. H., Tauberian theorems and the central limit theorem.  {\sl Ann. Probab.} {\bf 9} (1981), 221-231. \\
%  MR0606985 (82f:40010, H. Kesten).
\ni [Bin1984a] Bingham, N. H.,  On Euler and Borel summability.  {\sl J. London Math. Soc.} (2) {\bf 29} (1984), 141-146. \\
% [32] MR0734999 (85k:40010, E. Smet).
\ni [Bin1984b] Bingham, N. H., On Valiron and circle convergence.  {\sl Math. Z.} {\bf 186} (1984), 273-286. \\
% [33] MR0741307 (86g:40008, Amnon Jakimovski).
\ni [Bin1984c] Bingham, N. H., Tauberian theorems for summability methods of random-walk type.  {\sl J. London Math. Soc.} (2) {\bf 30} (1984), 281-287. \\
% [34] MR0771423 (86f:60085, E. Cs\'aki).
\ni [Bin1989] Bingham, N. H., Moving averages.  {\sl Almost Everywhere Convergence I}  (ed. G.A. Edgar \& L. Sucheston) 131-144, Academic Press, 1989. \\
% [46] MR1035241 (91c:60032, J. Steinebach).
\ni [Bin2008] Bingham, N. H., Tauberian theorems and large deviations.  {\sl Stochastics} {\bf 80} (Special Issue: A Festschrift for Priscilla Greenwood) (2008), 143-149. \\
\ni [Bin2019] Bingham, N. H.,  Riesz means and Beurling moving averages.  {\sl Risk and Stochastics: Ragnar Norberg} (Memorial Volume, ed. P. M. Barrieu) 159-172, World Scientific, 2019. \\
% [143] MR4404375;  arXiv:1502.07494. 
\ni [Bin2021] Bingham, N. H., The life, work and legacy of P. L. Chebyshev.  {\it P. L. Chebyshev - 200} (ed. A. N. Shiryaev).  {\sl Teor. Veroyatnost. Primen.} {\bf 66} (2021), 636-656. \\
% MR4331213; arXiv:2111.12551.
\ni [BinGa2015] Bingham, N. H. and Gashi, Bujar, Logarithmic moving averages. {\sl J. Math. Anal. Appl.} 421 (2015). \\
% [125]  MR3258350 (Kajal Khatri).  arXiv:1408.1301. 
\ni [BinGa2017] Bingham, N. H. and Gashi, Bujar, Voronoi means, moving averages and power series.  {\sl J. Math. Analysis and Applications} {\bf 449}.1 (2017), 682-696.  \\
% [132] MR3595227 (Kejal Khatri). arXiv:1607.02455.
\ni [BinGo1983] Bingham, N. H. and Goldie, C. M.,  On one-sided Tauberian conditions. {\sl Analysis} {\bf 3} (1983), 159-188. \\
% MR0756113 (85m:40004 Y. Sitaraman).
\ni [BinGo1988]  Bingham, N. H. and Goldie, C. M., Riesz means and self-neglecting functions.  {\sl Math. Z.} {\bf 199} (1988), 443-454.\\
\ni [BinGT1987]  Bingham, N. H., Goldie, C. M. and Teugels, J. L., {\sl Regular variation}.  Encycl. Math. Appl. {\bf 27}, Cambridge University Press, 1987. \\  
\ni [BinO2013] Bingham, N. H. and Ostaszewski, A. J., Steinhaus theory and regular variation: De Bruijn and after.  {\sl Indagationes Mathematicae} (N. G. de Bruijn Memorial Issue) {\bf 24} (2013), 679-692. \\
%BinO18  MR3124800 (D. J. Djur\u ci\'c).
\ni [BinO2014] Bingham, N. H. and Ostaszewski, A. J., Beurling slow and regular variation.  {\sl Transactions of the London Mathematical Society} {\bf 1} (2014), 29-56. \\
%BinO19
\ni [BinO2016] Bingham, N. H. and Ostaszewski, A. J., Beurling moving averages and approximate homomorphisms.  {\sl Indagationes Mathematicae} {\bf 27} (2016), 601-633. \\
\ni [BinO2020a] Bingham, N. H. and Ostaszewski, A. J., General regular variation, Popa groups and quantifier weakening. {\sl J. Math. Anal. Appl.} {\bf 483} (2020), 123610. \\
% BinO28  MR4029007 (Michael Voit). arXiv:1901.05996.
\ni [BinO2020b] Bingham, N. H. and Ostaszewski, A. J., Sequential regular variation: extensions to Kendall's theorem.  {\sl Quart. J. Math.} {\bf 71}(4) (2020), 1171-1200. \\ 
% BinO30 MR4186515 (Dragan \u Zimovir Djur\u ci\'c). arXiv:1901.07060.  
\ni [BinO2021] Bingham, N. H. and Ostaszewski, A. J., Extremes and regular variation.  Ch. 7, {\sl A lifetime of excursions through random walks and L\'evy processes}.  Progr. Prob. {\bf 78}, 121-137, Birkh\"auser, 2021 (A volume in honour of Ron Doney's 80th birthday, ed. L. Chaumont and E. A. Kyprianou). \\  
%BinO34 arXiv:2001.05420.
\ni [BinT1986] Bingham, N. H. and Tenenbaum, G.,  Riesz and Valiron means and fractional moments. {\sl Math. Proc. Cambridge Phil. Soc.} {\bf 99} (1986), 143-149. \\ 
% [37] MR0809509 (86m:40011, E. Omey).
\ni [Blo1976] Bloom, S., A characterisation of $B$-slowly varying functions.  {\sl Proc. Amer. Math. Soc.} {\bf 54} (1976), 243-250. \\
\ni [Bor1899] Borel, E., M\'emoire sur les s\'eries divergentes.  {\sl Ann. Sci. Ecole Norm. Sup.} {\bf 16} (1899), 9-136. \\
\ni [ChaM1952] Chandrasekharan, K. and Minakshisundaram, S., {\sl Typical means}.  Oxford University Press, 1952. \\
\ni [Cho1973] Chow, Y.-S., Delayed sums and Borel summability of independent, identically distributed random variables.  {\sl Bull. Inst. Math. Acad. Sinica} {\bf 1}, 207-220. \\
\ni [ChoL1973] Chow, Y.-S. and Lai, T.-L., Limiting behaviour of weighted sums of independent random variables.  {\sl Ann. Prob.} {\bf 1} (1973), 810-824. \\
\ni [Cra1938] Cram\'er, H., Sur un nouveau th\'eor\`eme-limite de la th\'eorie des probabilit\'es.  {\sl Act. Sci. et Ind.} {736} (1938), 5-32 (reprinted in {\sl Collected Works Vol. II}, 895-913, Springer, 1994). \\
\ni [DemZ1998] Dembo, A. and Zeitouni, O., {\sl Large deviations techniques and applications}, 2nd ed., Springer, 1998 (1st ed., Jones and Bartlett, 1993). \\
\ni [DenD] D\'eniel, D. and Derriennic, Y., Sur la convergence presque sure, au sens de Ces\`aro d'ordre $\a$, $0 < \a < 1$, des variables al\'eatoires ind\'ependantes et identiquement distribu\'ees.  {\sl Prob. Th. Rel. Fields} {\bf 79}, 629-636. \\
\ni [DeuS1989] Deuschel, J.-D. and Stroock, D. W., {\sl Large deviations}.  Academic Press, 1989. \\
\ni [DupE1997] Dupuis, P. and Ellis, R. S., {\sl A weak convergence approach to the theory of large deviations}.  Wiley, 1997. \\
\ni [ElA2017] El Arrouchi, M., Characterization of tail distributions based on record values by using Beurling's Tauberiqan theorem.  {\sl Extremes} {\bf 20} (2017), 111-120. \\
\ni [EmbK1997] Embrechts, P., Kl\"uppelberg, C. and Mikosch, T., {\sl Modelling extremal events for insurance and finance}.  Springer, 1997. \\
\ni [Fel1971] Feller, W., {\sl An introduction to probability theory and its applications, Volume II}, 2nd ed., Wiley, 1971. \\
\ni [Gai1953] Gaier, D., Der allgemeine L\"uckenumkehrsatz f\"ur das Borel-Verfahren.  {\sl Math. Z.} {\bf 88} (1965), 410-417. \\
\ni [Hal1992] Hall, P., {\sl The bootstrap and Edgeworth expansion}.  Springer, 1992. \\
\ni [Har1949] Hardy, G. H., {\sl Divergent series}.  Oxford University Press, 1949. \\
\ni [HarL1916] Hardy, G. H. and Littlewood, J. E., Theorems concerning the summability of series by Borel's exponential method.  {\sl Rend. Circ. Mat. Palermo} {\bf 41} (1916), 36-53 (reprinted in {\sl Collected Papers of G. H. Hardy VI}, 609-628, Oxford University Press, 1969). \\
\ni [Haw1980] Hawkes, J., Probabilistic behaviour of some lacunary series.  {\sl Z. Wahrschein. verw. Geb.} {\bf 53} (1980), 21-33. \\
\ni [IbrL1971] Ibragimov, I. A. and Linnik, Yu. V., {\sl Independent and stationary sequences of random variables}.  Wolters-Noordhoff, 1971. \\
\ni [Kor1959] Korevaar, J., Pansions and the theory of Fourier transforms.  {\sl Trans. Amer. Math. Soc,} {\bf 91} (1959), 53-101. \\
\ni [Kor1965] Korevaar, J., Distribution proof of Wiener's Tauberian theorem.  {\sl Proc. Amer. Math. soc.} {\bf 16} (1965), 353-355. \\
\ni [Kor2004] Korevaar, J., {\sl Tauberian theory, A century of developments}.  Grundl. math. Wiss. {\bf 329}, Springer, 2004. \\
\ni [Kor2013] Korevaar, J., Early work of N. G. (Dick) de Bruijn in analysis and some of my own.  {\sl Indag. Math.} {\bf 24} (2013), 668-678. \\
\ni [Kor2015] Korevaar, J., Johannes Gualtherus van der Corput (4 September 1890 - 13 September 1975).  {\sl Indag. Math.} {\bf 26} (2015), 715-722. \\
\ni [KorvA1949] Korevaar, J., van Aardenne-Ehrenfest, T. and de Bruijn, N. G., A note on slowly oscillating functions.  {\sl Nieuw Arch. Wiskunde} {\bf 23} (1949), 77-86. \\
\ni [Lai1974] Lai, T.-L., Summability methods for independent, identically distributed random variables.  {\sl Proc. Amer. Math. Soc.} {\bf 45} (1974), 253-261. \\
\ni [Lev1940] Levinson, N., {\sl Gap and density theorems}.  AMS Colloq. Publ. {\bf XXVI}, Amer. Math. Soc., 1940. \\
\ni [Lor1955] Lorentz, G. G., Borel and Banach properties of methods of summation.  {\sl Duke Math. J.} {\bf 22}, 129-141. \\
\ni [Mat1964] Matuszewska, W., On a generalization of regularly increasing functions.  {\sl Studia Math.} {\bf 24} (1964), 271-279. \\
\ni [MeyK1949] Meyer-K\"onig, W., Untersuchungen \"uber einige verwandte Limitierungsverfahren.  {\sl Math. Z.} {\bf 52} (1949), 257-304. \\
\ni [MeyKZ1960] Meyer-K\"onig, W. and Zeller, K., On Borel's method of summability.  {\sl Proc. Amer. Math. Soc.} {\bf 11} (1960), 307-314. \\
% 15.8.1984: Discrete Borel has no pure gap theorem. 
\ni [Ost2015] Ostaszewski, A. J., Beurling regular variation, Bloom dichotomy, and the Go{\l}ab-Schinzel functional equation.  {\sl Aequat. Math.} {\bf 89} (2015), 725-744. \\
\ni [Pet1975] Petrov, V. V., {\sl Sums of independent random variables}.  Ergeb. Math. {\bf 82}, Springer, 1975. \\ 
\ni [Tau1897] Tauber, A., Ein Satz aus der Theorie der unendlichen Reihen.  {\sl Monatsh. Math. u. Phys.} {\bf 8} (1897), 273-277. \\
\ni [Ten1980] Tenenbaum, G., Sur le proc\'ed\'e de sommation de Borel et la r\'epartition du nombre des facteurs premiers des entiers.  {\sl Enseignement Math.} {\bf 26} (1980), 225-245. \\
\ni [Tur1984] Tur\'an, P., {\sl On a new mthod of analysis and its applications}.  Wiley, 1984. \\
\ni [Wei1959] Weiss, M., The law of the iterated logarithm for lacunary trigonometric series.  {\sl Trans. Amer. Math. Soc.} {\bf 91} (1959), 444-469. \\
\ni [Wie1932] Wiener, N., Tauberian theorems.  {\sl Ann. Math.} {\bf 33} (1932), 1-100 (reprinted in {Generalized harmonic analysis and Tauberian theorem}, MIT Press, 1966). \\
\ni [ZelB1970] Zeller, K. and Beekmann, W., {\sl Theorie der Limitierungsverfahren}, 2nd ed.  Springer, 1970. \\

\ni N. H. Bingham, Mathematics Department, Imperial College, London SW7 2AZ; n.bingham@ic.ac.uk \\ 

\end{document}